\input amstex
\input psfig.sty
\documentstyle{amsppt}
\magnification=1200
\baselineskip=16pt
\lineskip=4pt
\lineskiplimit=4pt
\hsize=5.5in
\vsize=7.0in
\parskip=5pt
\TagsOnRight
\define\hcenterbox#1{$$\vbox{#1}$$}
\define\larrow{\leftarrow}

\define\mck{M\raise0.6ex\hbox{\kern-1pt{c}}Kay }
\define\MCK{M\raise0.6ex\hbox{\kern-1pt{c}}KAY}

\define\Cox{\mathop{{\text{Cox}}}\nolimits}
\define\bA{{\bf{A}}}
\define\bD{{\bf{D}}}
\define\bE{{\bf{E}}}
\topmatter
\title Semi-affine Coxeter-Dynkin graphs and $G \subseteq SU_2(C)$
\endtitle
\author 
JOHN \MCK
\endauthor
\affil Concordia University, \ Montr\'eal, Qu\'ebec, \ Canada \ H3G 1M8
\endaffil
\footnote""{Partially supported by NSERC and FCAR grants.}
\abstract\nofrills{{\bf Abstract.} 
The semi-affine Coxeter-Dynkin graph is introduced, generalizing both
the affine and the finite types.}
\endabstract

\endtopmatter

\document

\subheading{Semi-affine graphs}

It is profitable to treat the so-called Coxeter-Dynkin diagrams as graphs.
A classification of finite graphs with an adjacency matrix having 2 as the
largest eigenvalue is made 
in a paper of John Smith \cite {JHS}. It is in a combinatorial context
and no reference is made to Coxeter-Dynkin diagrams there. This maximal
eigenvalue property is a defining property of the affine diagrams. What 
is introduced in this note is a more weakly constrained graph, and we 
examine its eigenvalues
and interpret the rational functions which arise in terms of my 
correspondence \cite {M1,K}. Since these semi-affine graphs do not have
symmetrizable matrices, this appears to imply a connection with
singularities rather than Lie algebras. 

Here we shall deal only with those of type \bA, \bD, and \bE. Undirected
edges are treated as a pair of edges directed in opposing directions as in
\cite{FM,M1,M2}. By so doing, we can introduce the semi-affine graph
which may be defined in terms of a graph of finite type with an additional
edge (two for \bA-type) directed toward the affine node; equivalently it may
be defined as an affine graph with any undirected edge connecting the affine
node replaced by an directed edge directed toward the affine node. This is 
done by removing one of the two opposed directed edges.

Effectively the semi-affine graph generalizes both the affine and finite
type graph since the affine node acts as a sink, and when
weighted at the nodes, it satisfies the same constraints as the affine
graph except the constraint imposed by the additional directed edge(s)
in the affine graph. Note that weighting the extra node with zero yields the 
same constraints as the finite type graph.
In some sense the semi-affine graph lies intermediate between the finite and
affine graphs yet generalizes both.

It is helpful to have a simple example at hand: for this we choose the
type \bD$_4$. The finite, semiaffine, and affine graphs of this type are:

\hbox{\qquad}\psfig{file=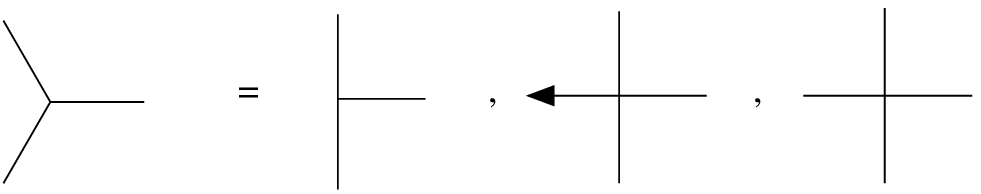,width=10cm,height=2cm}



To each of the $r+1 $ nodes, $i$, we attach a weight, $n_i$, satisfying

$$
     t n_i  = \sum_{j \larrow i}  n_j,                           \tag{$*$}
$$
 summed over the successor nodes, $j$, of $i$. Consistent with the 
geometrical interpretation of summing
over successor nodes, see \cite{M1}, we note that $n_0$ occurs in the 
right side only. We initially normalize so that the affine node, $n_0$, is given
the weight, 1. This yields weights of $t/(t^2-3)$ for the central
node, and $1/(t^2-3)$ for the other nodes. We now renormalize to make the
weights polynomials in $t$ by multiplying by the denominator, $t^2-3$.

Now $t^2-3$ is the minimal polynomial of $2\cos(2\pi/2h)$ where $h$ is the
Coxeter number (= $6$ here). This suggests writing $t = q+1/q$ and clearing
denominators. This yields respective weights of $1-q^2+q^4, q+q^3 $, and
$q^2 $ for the affine, central, and other nodes.  These are the numerators of 
generalized Molien series (to be described below) when written as a rational 
function of $q$ with no common factor.

Normalizing the affine node to $1+q^h = 1+q^6$ by multiplying by $1+q^2$,
we find that the weights: $1+q^6,q+2q^3+q^5$, and three of $q^2+q^4$, are 
the numerators of the generalized Molien series, \cite{NJAS}, but now written
with standard denominator for finite $G \subset SU_2(C)$, 

More generally, we find each $n_i$ is a 
rational function of $t$. For a finite graph we have a common denominator 
polynomial, $\Cox(t)$, which is the minimal polynomial of $2\cos(2\pi/2h)$,
for $h$, the Coxeter number of the Lie algebra associated with the graph.
We put $t = q + q^{-1}$ and clear denominators
to obtain polynomials, $n_i(q)$, for weights normalized
so that $n_0(q) = 1 + q^h$. 

To a representation $R$, and an irreducible character $\chi_i$, of $G$, 
the generalized Molien series is defined by:
$$
	m_i(G) 
	= \frac{1}{|G|} \sum_{x \in G} \frac{\chi_i(x)}{\det(I-R(x)q)},
$$
which, see \cite{NJAS}, for finite $R(G) \subset SU_2(C)$, 
can be expressed in standard form as
$$
	m_i(G) = \frac{N_i(q)}{(1-q^a)(1-q^b)}
$$
with  
$ab = 2|G|$, $a+b = h+2$, $h$ the Coxeter number assigned
to $G$ by the \mck correspondence, and $\{\chi_i\}$ 
($\chi_0(x) = 1$),
being the
set of irreducible characters of $G$.
We find that
$$
	(q+q^{-1})m_i(G) = \sum_{j \larrow i} m_j(G)
$$
and $N_i(q) = n_i(q)$. 

The semi-affine graph has characteristic polynomial $t^d\Cox(t)$
of total degree = rank + 1,
where degree $\Cox(t) = \varphi (2h)/2$ and $\varphi$ is Euler's function.

\subheading{Specialization}

The condition $N_0(q) = 1 + q^h = 0$ yields numeric weights
for the finite type graphs.

We may instead impose the extra condition obtained by making the semi-affine
graph into an affine one --- this gives numeric weights for the affine graph,
and these values are those for which the denominator of the standard form
of the Molien series vanishes.

\hcenterbox{
\halign{ 
# \hfil\ : \quad & \hfil $#$ &= $#$ \hfil 
	& \quad implying \quad \hfil $#$ &= $#$ \hfil \cr
\bA & (q+q^{-1})(1+q^h) & 2(q+q^{h-1})          &     (1-q^2)(1-q^h) & 0;\cr
\bD & (q+q^{-1})(1+q^h) & q+q^3+q^{h-3}+q^{h-1} & (1-q^4)(1-q^{h-2}) & 0;\cr
\bE & (q+q^{-1})(1+q^h) & q+q^{a-1}+q^{b-1}+q^{h-1} & (1-q^a)(1-q^b) & 0.\cr
}
}

It is useful to note:

1. Chains starting at the affine node have the smallest $q$ exponent increasing
   by 1 at each successive node - similarly the largest exponent decreases by 1.

2. For even $h$, the $q$ exponents are alternately all even and all odd at 
   adjacent nodes. For odd $h$ (=$A_m$, m even), each odd exponent pairs with an
   even one.

3. The number of powers of $q$ at a node is half the number of powers of $q$ 
   summed over adjacent nodes.

For the sake of brevity, I give the $q$ exponents along the longest chain 
starting at the affine node. Other nodes are either given last, or are 
determined by a graph symmetry fixing the affine node.

$A_m: n_k=q^k+q^{h-k}, k=0,...,m, h=m+1$.

$D_m$: At the tips, $n_0 = 1+q^h, n_1=q^2+q^{h-2}$ and $n_m=n_{m+1} = q^{m-2}+q^m$.

      The $m-3$ central nodes are weighed  $q^k+q^{h-k}+q^{k+2}+q^{h-k+2}, k=
      1,...,m-3, h=2m-2$.

$E_8: (0+30),(1+11+19+29),(2+10+12+18+20+28),(3+9+11+13+17+19+21+27),
      (4+8+10+12+14+16+18+20+22+26),(5+7+9+11+13+2\times 15+17+19+21+23+25),
      (6+8+12+14+16+18+22+24),(7+13+17+23) + (6+10+14+16+18+20+24)$.

$E_7: (0+18),(1+7+11+17),(2+6+8+10+16),(3+5+7+2\times 9+11+13+15),(4+6+8+10+12+14),
      (5+7+11+13),(6+12) + (4+8+10+14)$.

$E_6: (0+12),(1+5+7+11),(2+4+2\times 6+8+10),(3+5+7+9),(4+8)$ + symmetry.

\subheading{Interpretation of $q$}

The polynomials 
$n_i(q)$ are self-reciprocal since we start with rational functions of 
$t = q + q^{-1}$. This may be interpreted in terms of Poincar\'e duality.

The coefficients of $q^k$ in the numerators $N_i$ of $m_i(G)$ count the 
dimensions of certain fixed spaces under the group action \cite {NJAS};
this exhibits $q$ as a dimension-shifter. One may also interpret the main
equation ($*$) as a trace equation in which we see $q+q^{-1}$ as the
trace of an element in $SU_2(C).$


\enddocument